\documentclass[10pt]{amsart}
\usepackage[all]{xy}
\usepackage[english]{babel}

\begin{document}

\title{Inheriting of chaos in nonautonomous  dynamical systems}

\author{M. \v Stef\'ankov\'a}

\address{M. \v Stef\'ankov\'a, Mathematical Institute, Silesian University, 746
01 Opava, Czech Republic}

\email{marta.stefankova@math.slu.cz}

\thanks{ The research was supported, in part, by the grant  GA201/10/0887 from the Czech Science Foundation.  Support of this institution is gratefully acknowledged.}

\begin{abstract}  We consider nonautonomous discrete dynamical systems $\{ f_n\}_{n\ge 1}$, where every $f_n$ is a surjective continuous map $[0,1]\to [0,1]$ such that $f_n$ converges uniformly to a map $f$. We show, among others, that if $f$ is chaotic in the sense of Li and Yorke then the nonautonomous  system $ \{ f_n\}_{n\ge 1}$  is Li-Yorke chaotic as well, and that the same is true for  distributional chaos. If $f$ has zero topological entropy then the nonautonomous system inherits its infinite $\omega$-limit sets.
\\[.3cm]
\end{abstract}

\maketitle

\pagestyle{myheadings}
\markboth{{\sc M. \v Stef\'ankov\'a}}{{\sc Nonautonomous dynamical systems}}

\section{Introduction}

Let $(X,\rho )$ be a compact metric space, and $I=[0,1]$. Denote by $\mathcal C(X)$ the class of continuous maps $X\to X$, and let $\mathcal C$ stand for $\mathcal C(I)$. By $\mathbb N$ and $\mathbb N_0$ we denote the set of positive, or nonnegative integers, respectively.  $(X,f)$, with $f\in \mathcal C(X)$, is  a {\it topological dynamical system\/}. A {\it nonautonomous (discrete dynamical) system} is a pair $(X,\{ f_n\}_{n\ge 1})$, where $f_n\in\mathcal C (X)$, $n\in\mathbb N$; following \cite{KS} we denote this system by $(X, f_{1,\infty})$. The {\it trajectory} of an $x\in X$ in this system is the sequence $\{x_n\}_{n\ge 0}$, where $x_0=x$ and $x_n=(f_n\circ f_{n-1}\circ\cdots\circ f_1)(x)$. The set of limit points of the trajectory of a point $x$ is its {\it $\omega$-limit set}; we denote it by $\omega_{f_{1,\infty}}(x)$. If $f_n=f$ for every $f_n\in f_{1,\infty}$ then  $(X,f_{1,\infty})=(X,f)$.

Nonautonomous  systems are closely related  to { skew-product maps} $F: X\times Y\to X\times Y$, with $X,Y$ compact metric spaces; for details see, e.g., \cite{KS},  a pioneering work dealing with nonautonomous systems, motivated just by open problems concerning skew-product maps. In particular, \cite{KS}  deals with  topological entropy which can be for nonautonomous systems defined similarly as for the autonomous ones. We denote by $h(f)$ or $h(f_{1,\infty})$ the topological entropy of a map $f$, or $f_{1,\infty}$, respectively.

Let $\{x_n\}_{n\ge 0}$, $\{y_n\}_{n\ge 0}$ be trajectories of points $x,y\in X$, and $\varepsilon >0$. Then $(x,y)$ is an $\varepsilon$-{\it Li-Yorke pair} if
$\limsup _{n\to\infty} \rho (x_n,y_n)\ge\varepsilon$ and $ \liminf _{n\to\infty}\rho (x_n, y_n)=0$. For $x,y\in X$ define $\Phi_{xy}, \Phi^*_{xy}:(0,\infty) \to I$  by  
\begin{equation}
\label{equ12}
\Phi _{xy}(t):=\liminf _{n\to\infty}\frac 1n\#\{ 0\le j<n; \rho (x_j,y_j)<t\}, \ \Phi ^*_{xy}(t):=\limsup _{n\to\infty}\frac 1n\#\{ 0\le j<n; \rho (x_j,y_j)<t\}. 
\end{equation}
The system  $(X,f)$ or $(X, f_{1,\infty})$ is {\it Li-Yorke chaotic}, briefly LYC, if there is an $\varepsilon >0$, and an uncountable {\it scrambled set} $S$ such that every distinct points  $x,y\in S$ form an $\varepsilon$-Li-Yorke pair; it is {\it distributionally chaotic}, briefly DC1, if there there is an $\varepsilon>0$, and an uncountable set $S$ such that for every distinct points  $x,y\in S$, $\Phi _{xy}(\varepsilon)=0$ and $\Phi^*_{xy}\equiv 1$. Notice that $(I,f)$ is  DC1 if and only if $h(f)>0$, see  \cite{ScSm}.
\medskip

Our paper is inspired by \cite{KS} (see also \cite{C}) where relations between systems $(I,f_{1,\infty})$ and $(I,f)$ such that $f_{1,\infty}$ uniformly converges to $f$ are considered. Since a single constant function in $f_{1,\infty}$ can destroy more complex behavior, even in the case when the limit system $(I,f)$ has complicated dynamics, {\it in this paper we  assume that $f$ and all maps in $f_{1,\infty}$ are surjective}. With this condition, for example, $h(f)>0$ implies $h(f_{1,\infty})>0$  \cite{C},  without it we have only $h(f)\ge h(f_{1,\infty})$ \cite{KS}.  Consequently, if $h(f)>0$ then it is possible to show directly that the nonautonomous system is DC1 (we obtain this result indirectly from  Theorem  3.2). Therefore our paper is devoted to systems with zero topological entropy. The proofs are based on  \lq\lq classical\rq\rq  papers concerning chaos and structure of $\omega$-limit sets of maps  $f\in\mathcal C$ with $h(f)=0$,  \cite{Sh}, \cite{S}, \cite{FShSm}, \cite{BS}. Our main result is Theorem C; we show  that $(I,f_{1,\infty})$ is LYC if $(I,f)$ is LYC. In some cases, the nonautonomous system inherits stronger forms of chaos (Theorem B) and infinite $\omega$-limit sets (Theorem A). Note that $(I,f)$ need not be LYC  or DC1 if $(I,f_{1,\infty})$ is, see \cite{FPS}. Theorem A is interesting in itself and makes possible to prove other results more transparently.\\

{\bf Theorem A.} {\it  Let $(I, f_{1,\infty})$ be a surjective nonautonomous system, and let $f_{1,\infty}$ converge uniformly to a map $f$. If $h(f)=0$ then every infinite $\omega$-limit set of $f$ is an $\omega$-limit set for $f_{1,\infty}$.}\\

{\bf Theorem B.} {\it Let $(I, f_{1,\infty})$ be a surjective nonautonomous system, and $f_{1,\infty}$ converge uniformly to a map $f$. Then $(I,f_{1,\infty})$ is DC1 if one of the following conditions is satisfied:

(i) $h(f)>0$ (or equivalently,  $f$ is DC1); 

(ii) $f$  has a minimal set $\widetilde\omega$ such that $f|_{\widetilde\omega}$ is  not Lyapunov stable.}
\smallskip

\noindent Recall that $f$ is {\it Lyapunov stable} if for every $\varepsilon >0$ there is a $\delta >0$ such that $|x-y|<\delta$ implies $|f^n(x)-f^n(y)|<\varepsilon$, for every $n$.\\

{\bf Theorem C.} {\it Let $(I, f_{1,\infty})$ be a surjective nonautonomous system, and $f_{1,\infty}$ converge uniformly to a map $f$. If $f$  is LYC  then also $(I,f_{1,\infty})$ is LYC.}\\

{\bf Remarks.} Obviously, Theorem A is not valid for finite $\omega$-limit sets.  Theorems B and C cannot be strengthened in the sense that  the non-autonomous system inherits chaos with extremal properties like big scrambled sets. For example, a map in $\mathcal C$ can have DC1 scrambled set with complement of zero Hausdorff dimension \cite{OS}, but this need not be inherited by a nonautonomous system, see  \cite{Dv}. Theorem B is interesting since there are functions $f\in\mathcal C$ with $h(f)=0$ satisfying condition (ii), see \cite{FShSm} or \cite{BS} .
\medskip

\section{Proof of Theorem A}

A  compact set $A\subseteq X$ is {\it $f$-periodic  of period} $m$, where $f\in\mathcal C(X)$, if $f^j(A)$ are pairwise disjoint, for $0\le j<m$,  and $f^m(A)=A$.
\medskip
 
  {\bf Theorem 2.1.} (See \cite{S}.) {\it Let $f\in\mathcal C$ with $h(f)=0$,  and let $\widetilde\omega$ be an infinite $\omega$-limit set of $f$.
Then there is a system $\{J(k,n); 0\le k<2^n\}_{n\ge 0}$ of $f$-periodic intervals in $I$ such that, for any $k, n\in\mathbb N_0$, 
\smallskip

(i) $f(J(k,n))=J(k+1,n)$ where $k+1$ is taken ${\rm mod} \ 2^n$;
\smallskip
  
(ii) $J(k,n)$ has period $2^n$;
\smallskip

(iii) $J(k,n+1)\cup J(2^n+k, n+1)\subset J(k,n)$;
\smallskip

 (iv) $\widetilde\omega\subset\bigcup _{0\le k<2^n} J(k, n) =: O_n$.}\\

\noindent Obviously we may assume that the intervals $J(k,n)$ are the minimal ones in the sense of inclusion. In this case, the collection of all $J(k,n)$ is the {\it system associated to} $\widetilde\omega$; we denote it by $\mathcal J_f(\widetilde\omega)$, or simply by $\mathcal J$. 
The  system
\begin{equation}
\label{sup01}
 \{ \widetilde\omega (k,n):=J(k,n)\cap\widetilde\omega\}_{0\le k<2^n}, \  k,n\in\mathbb N, 
 \end{equation}
 is {\it the cyclic decomposition of} $\widetilde\omega$ {\it of degree} $n$. Since $f(\widetilde\omega)=\widetilde\omega$, by Theorem 2.1 
 \begin{equation}
 \label{sup02}
\{\widetilde\omega(k,n)\}_{0\le k<2^n} \ \text{forms an} \ f \text{-periodic orbit} \ \text{of period} \   2^n,  \ \text{and} \  \bigcup _{0\le k<2^n} \widetilde\omega(k,n)=\widetilde\omega.
\end{equation}
For the cyclic decomposition (\ref{sup02}), and $k,n\in\mathbb N$, $0\le k<2^n$, let $K(k,n)\subset J(k,n)$ be the compact interval between the sets $\widetilde\omega(k,n+1)$ and $\widetilde\omega(2^n+k,n+1)$ (which are neighbor sets in the cyclic decomposition of $\widetilde\omega$ of degree $n+1$). This $K(k,n)$ is {\it a complementary interval to} $\widetilde\omega$ of degree $n$.\\

{\bf Lemma 2.2.} {\it Assume that $f\in\mathcal C$, $h(f)=0$, and  $\widetilde\omega$ is infinite $\omega$-limit set of $f$.  Then
\begin{equation}
\label{sup03}
f^{2^{n+2}}(K(k,n))\supset\widetilde\omega (k,n), \ k, n\in\mathbb N, \ 0\le k<2^n.
\end{equation}}
\smallskip
{\bf Proof.} We may assume that $\widetilde\omega(k,n+1)<\widetilde\omega(2^n+k,n+1)$ where $<$ indicates the natural ordering of (disjoint) sets. Let $\widetilde\omega(k_0,n+2)< \widetilde\omega(k_1,n+2)<\widetilde\omega(k_2,n+2)<\widetilde\omega(k_3,n+2)$ be the sets from the cyclic decomposition of $\widetilde\omega$ of degree $n+2$ contained in $\widetilde\omega(k,n)$ so that $\widetilde\omega(k_0,n+2)\cup\widetilde\omega(k_1,n+2)\subset \widetilde\omega(k,n+1)$ and $\widetilde\omega(k_2,n+2)\cup\widetilde\omega(k_3,n+2)\subset \widetilde\omega(2^n+k,n+1)$. Since $\widetilde\omega(k,n)$ is an $\omega$-limit set of $f^{2^n}$ and no point in $\widetilde\omega$ is periodic, the interval between $f^{2^n}(u)$ and $f^{2^n}(v)$, where $u,v$ are the endpoints of $K(k,n)$,  must contain one of the sets $\widetilde\omega(k_j,n+2)$, $0\le j<4$. Consequently, the interval between $f^{2^{n+1}}(u)$ and $f^{2^{n+1}}(v)$ contains two of the sets, the interval between $f^{3\cdot 2^{n}}(u)$ and $f^{3\cdot2^{n}}(v)$ three of the sets, and (\ref{sup03}) follows.
$\hfill\Box$\\

{\bf Lemma 2.3. } (Itinerary lemma.)  {\it Let $f_{1,\infty}$ be a sequence of maps in $\mathcal C(X)$, and  $F_{1,\infty}$ a sequence of nonempty compact subsets of $X$ such that, for every $n\in\mathbb N$, $f_n(F_n)\supseteq F_{n+1}$. Then there is an $x$ such that $x_n\in F_n$, $n\in N$, where $x_1,x_2,\cdots$ is the trajectory of $x$ in the nonautonomous system.
\smallskip}

{\bf Proof} is easy. $\hfill\Box$\\

{\bf Lemma 2.4.} {\it Assume $f\in\mathcal C$ with $h(f)=0$,  $\omega _f(z)=:\widetilde\omega$ is infinite, and  $p$ is an isolated point of $\widetilde\omega$. Then there is a cluster point  $a_p$ of $\widetilde\omega$ such that the interval $J_p$ with endpoints $p$ and $a_p$ is a wandering interval (i.e., $f^i(J_p)\cap f^j(J_p)=\emptyset$ if $i\ne j$), and for every neighborhood $U$ of $p$ and every $m\in\mathbb N$ there is a $q\in\mathbb N$ divisible by $2^m$ such that  $f^q(U)$ is a neighborhood of $J_p$.}\\

{\bf Proof.} This result, in a different setting, is a part of Lemma 2.9 in \cite{S}. For convenience, we provide an outline of the argument. Let $J(k_0,0)\supset J(k_1,1)\supset \cdots\supset J(k_n,n)\supset\cdots$ be the intervals in $\mathcal J_f(\widetilde\omega)$ containing $p$. Then $\bigcap _{n\ge 0} J(k_n,n)=:J_p$ is a wandering interval with endpoints $p$ and $a_p\in\widetilde\omega$; moreover, $a_p$ is a cluster point of $\widetilde\omega$, see \cite{Sh} (cf. also \cite{BS}) so that $p$ is an endpoint of every $\widetilde\omega(k_n, n)$ with $n\ge n_0$. Let $z_{s(1)}, z_{s(2)}, z_{s(3)}, \cdots$ be a monotone subsequence of points in the trajectory of $z$ with $\lim _{i\to\infty}z_{s(i)}=p$; obviously, $z_{s(i)}\notin J_p$, $i\ge 1$. Let $z_{s(k)}\in U$. If $p$ is an endpoint of $\widetilde\omega (k_0,0)$ then, since $J_p$ is a wandering interval, $J_p$ is contained in the open interval $U^\prime$ with endpoints $z_{s(k+1)}$ and $f^{s(k+1)-s(k)}(p)$. To finish we may assume  $m\ge n_0$. Then $p$ is an endpoint of $\widetilde\omega (k_m,m)$, and  application of  the above process to $g:=f^{2^m}$ completes the argument. 
$\hfill\Box$\\

{\bf Proof of Theorem A.}  Denote by $P$ the set of isolated points of $\widetilde\omega$ and consider two possible cases.
 
\noindent {\it Case 1.} $P=\emptyset$ so that $\widetilde\omega$ is a minimal set of $f$.  For every $m,j\in\mathbb N$, $m\ge 1$, denote $f_{m}^j:=f_{m+j-1}\circ f_{m+j-2}\circ \cdots\circ f_{m+1}\circ f_m$. Since $f_{1,\infty}$ converges uniformly to $f$, by (\ref{sup03})  there is an $m(n)$ such that
\begin{equation}
\label{sup04}
f_m^{2^{n+2}}(K(k,n))\supset K(k,n+1)\cup K(2^n+k,n+1),  \ 0\le k<2^n, \ m\ge m(n),  \ k, m,n\in\mathbb N;
\end{equation} 
notice that $f_m^{2^{n+2}}(K(k,n))$ is a neighborhood of $K(k,n)$. Choose $c_n$ such that
\begin{equation}
\label{301}
m(n+1)-m(n)\le 2^nc_n, \ n, c_n\in\mathbb N, 
\end{equation}
where $m(n)$ is as in (\ref{sup04}). To simplify the notation let $K_n$ be the finite sequence $K(0,n),K(1,n),\cdots , K(2^n-1,n)$ of all $2^n$ intervals $K(k,n)$ of degree $n$.  We wish to apply Itinerary lemma to the sequence
\begin{equation}
\label{sup05}
F_{m(0),\infty}=\underbrace{K_0,K_0, \cdots , K_0}_\text{$c_0$-times},\underbrace{K_1,K_1, \cdots , K_1}_\text{$c_1$-times},\cdots ,\underbrace{K_n,K_n, \cdots , K_n}_\text{$c_n$-times},\cdots.
\end{equation}
Obviously, $f_j(F_j)\supseteq F_{j+1}$  if $f=f_j$ and,  by (\ref{301}),  if $F_j=K(k,n)$, for some $k, n$.  However, if the numbers $c_n$ are rapidly increasing, the inclusions will be satisfied \lq\lq approximately\rq\rq so that,  for every $j$, $F_{j+1}$ is contained in the  $\delta_j$ neighborhood of $f_j(F_j)$, where $\delta _j\to 0$. Apply Itinerary lemma to (\ref{sup05}), and  $(I,f)$ or $(I,f_{m(0),\infty})$, respectively, to get points $x$ and $x^\prime$ in $K(0,0)$. The trajectory of $x$ passes the sets in (\ref{sup05}) exactly, while the trajectory of $x^\prime$ hits exactly the sets $K(0,n)$.  The trajectories $\{x_j\}_{j\ge m(0)}$ and $\{x_j^\prime\}_{j\ge m(0)}$ of these points are proximal since $\delta _j\to 0$ so that  both must have the same $\omega$-limit set $\widetilde\omega^\prime$. But $\omega _f(x)=\widetilde\omega$ since by (\ref{sup05}) the trajectory of $x$ can have only finitely many members in the set $\bigcup _{0\le k<2^n}K(k,n)$ so that, by  Lemma 2.4,  $\omega_f(x)$ contains no isolated points. Since every $f_n$ is surjective,  $\omega_{f_{1,\infty}}(z)=\widetilde\omega$ for some $z\in I$.
\smallskip

{\it Case 2.} $P\ne\emptyset$. In the proof we need facts which are contained implicitly in the literature, see \cite{Sh}, \cite{S}, \cite{FShSm}, \cite{BS}; to make the proof self-contained, we recall some of them with brief arguments.  Let $\widetilde\omega=\omega _f(z)$, and let $\{z_j\}_{j\ge 0}$ be the trajectory of $z$. Since a point in $P$ cannot be periodic it has a preimage in $P$ so that  $P$ is countably infinite. Since the intervals in $\mathcal J$ are periodic, there are $j_n\in\mathbb N$ such that
\begin{equation}
\label{sup06}
z_j\in O_n \ \ \text{if} \ \ j\ge j_n, \ \text{and} \ z_j\notin O_n\setminus O_{n+1} \ \text {if} \ j\ge  j_{n+1},  \ \text{where} \ j_{n+1}>j_n, \ j,n\in\mathbb N, 
\end{equation}
where $O_n$ is the corresponding orbit of the intervals $J(k,n)$, $k\le 2^n-1$,  as in Theorem 2.1. To see this note that, by Theorem 2.1 (iv) and Lemma 2.4,  $\bigcap _{n \ge 1}O_n\setminus \widetilde\omega$ is the union of wandering intervals. It follows that for every $j$ there is a point $p_j$ such that the interval  with endpoints $p_j$ and $z_j$ intersects $\widetilde\omega$ exactly at $p_j$; denote this interval by $L_j$ and notice that $p_j$ need not be in $P$, since the image of an isolated point need not be isolated, see also \cite{BS}.  Obviously, $L_{j+1}$ has endpoints $z_{j+1}$ and $p_{j+1}:=f(p_j)$ so that $f(L_j)\supseteq L_{j+1}$. Since $L_j$ has just one point, $p_j$, in common with the wandering interval $J_{p_j}$, and $L_j\cup J_{p_j}$ is a neighborhood of $p_j$,  Lemma 2.4 applies to $U:=L_j$. Therefore 
\begin{equation}
\label{sup07}
f^{j_{n+1}-j_n}(L_{j_{n}})\supset K(k_n, n+1)\supset L_{j_{n+1}}. 
\end{equation}
For simplicity, denote by $\widetilde K_n$ the finite sequence $K(k_n,n+1), K(k_n+1,n+1), K(k_n+2,n+1),\cdots , K(2^{n+1}+k_n-1,n+1)$ which consists of the first $2^{n+1}$  sets in the $f$-trajectory of $K(k_n,n+1)$, and by $\widetilde L_n$ the finite sequence $L_{j_n}, L_{j_n+1},\cdots ,L_{j_{n+1}-1}$ of $j_{n+1}-j_n$ members of the $f$-trajectory of $L_{j_n}$. By Lemma 2.2 and (\ref{sup07}),  Itinerary lemma applied to $f$ and
\begin{equation}
\label{sup08}
\widetilde L_0, \underbrace {\widetilde K_0,\widetilde K_0,\cdots , \widetilde K_0}_\text{$c_0$-times}, \widetilde L_{1},\underbrace {\widetilde K_1,\widetilde K_1,\cdots , \widetilde K_1}_\text{$c_1$-times},\widetilde L_2,\cdots , \widetilde L_{n},\underbrace {\widetilde K_n,\widetilde K_n,\cdots , \widetilde K_n}_\text{$c_n$-times},\widetilde L_{n+1}, \cdots
\end{equation}
yields a point $x$ such that $\omega _f(x)=\widetilde\omega$ since  its trajectory passes through  $\widetilde L_0,\widetilde L_1, \cdots , \widetilde L_n, \widetilde L_{n+1},\cdots$. The inserted blocks  $\widetilde K_n$ in  (\ref{sup08}) contain only finitely many sets of type $K(i,n)$ which by Lemma 2.4  cannot generate new isolated points.  Similarly as in Case 1, replace the sequence $c_0,c_0,\cdots$ in (\ref{sup08}) by a more rapidly increasing sequence $\widetilde c_0,\widetilde c_1,\cdots$ if necessary, and apply Itinerary lemma to $f_{m,\infty}$ where $m$ is sufficiently large. This  gives a point $x^\prime$ such that $\omega_{f_{m,\infty}} (x^\prime)=\widetilde\omega$.
$\hfill\Box$

\section{Proofs of Theorems B and C.} 

 Recall that a map $f\in\mathcal C(X)$ has a {\it horseshoe} if there are disjoint nonempty compact sets $U,V$, and $m\in\mathbb N$ such that $f^m(U)\cap f^m(V)\supseteq U\cup V$. The following is a strictly weaker notion.

\medskip  

{\bf Definition 3.1.}  A map $f\in\mathcal C(X)$  has a {\it quasi horseshoe} if there are $\varepsilon >0$, compact sets $U_k, V_k$, and  positive integers $m_k$, for $k\in\mathbb N_0$, with the following properties:
\smallskip

(i) dist $(U_k, V_k)\ge \varepsilon$;
\smallskip

(ii) $\lim_{k\to\infty} {\rm diam}(U_k)=\lim_{k\to\infty} {\rm diam}(V_k)=0$;
\smallskip

(iii) $f^{m_k}(U_k)$ is a neighborhood of $U_k\cup U_{k+1}\cup V_{k+1}$,  
and $f^{m_k}(V_k)$ a neighborhood of $V_k\cup V_{k+1}\cup U_{k+1}$.\\

  {\bf Theorem 3.2.}  {\it Let $f, f_k\in\mathcal C(X)$ be surjective maps, for $k\in\mathbb N$, and let $f_{1,\infty}$ converge uniformly to $ f$.  If  $f$ has  a quasi horseshoe then $(X,f_{1,\infty})$ is distributionally (DC1) chaotic.}\\
  
  {\bf Proof.}  Keep the notation from Definition 3.1 and denote by $\widetilde U_k$ the finite sequence $U_k, f(U_k), f^2(U_k),$ $\cdots , f^{m_k-1}(U_k)$ of $m_k$ compact sets, and similarly with $\widetilde V_k$. Let $\Sigma _2=\{ 0,1\}^\mathbb N$. For $\alpha=\{a_k\}_{k\ge 0}\in\Sigma_2$ consider the itinerary
  \begin{equation}
  \label{401}
  I_\alpha: =\underbrace {B_0,  B_0,\cdots,  B_0}_{c_0{\text  -times}},\underbrace { B_1,  B_1\cdots,  B_1}_{c_1{\text  -times}}, \cdots, 
\underbrace { B_k,  B_k\cdots,  B_k}_{c_k{\text  -times}}, \cdots ,
 \end{equation}
 where
 \begin{equation}
 \label{402}
  B_k=\widetilde U_k \ \text {if} \ a_k=0, \ \text{and}  \  B_k=\widetilde V_k \ \text {if} \ a_k=1,  \ k\in\mathbb N_0.  
 \end{equation}
 If the numbers $c_k$ are sufficiently large then by Itinerary lemma, similarly as in the proof of Theorem A,  there is an  $x_\alpha\in U_0\cup V_0$  with itinerary $I_\alpha$ in $f_{1,\infty}$. Let $\Sigma^\prime_2\subset\Sigma _2$ be an uncountable set such that, for every distinct  $\{a_k\}_{k\ge 0}$ and $\{b_k\}_{k\ge 0}$ in $\Sigma^\prime _2$, we have $a_k=b_k$ for infinitely many $k$, and $a_k\ne b_k$ for infinitely many $k$; such a set exists, see, e.g., \cite{S}. Let $S=\{x_\alpha; \alpha\in \Sigma _2^\prime\}$ and assume that the numbers $c_k$ are  increasing so rapidly that $\lim _{k\to\infty} c_k/c_{k+1}=0$. Then it is easy to verify that $S$ is a DC1 scrambled set for $f_{1,\infty}$ such that, for every $x\ne y$ in $S$, $\Phi_{xy}(\varepsilon )=0$ and $\Phi^*_{xy}\equiv 1$.
  $\hfill\Box$
  
\medskip

The next theorem improves a result from  \cite{S} that a LYC map $f\in\mathcal C$ has similar system of intervals as in Definition 3.1 except that condition (iii) is replaced by $f^{m_k}(U_k)\cap f^{m_k}(V_k)\supset U_{k+1}\cup V_{k+1}$. The stronger property is necessary in our proof of Theorem 3.2.\\
 
{\bf Theorem 3.3.} {\it Let $f\in\mathcal C$ have a minimal set $\widetilde\omega$ such that $f|_{\widetilde\omega}$ is not Lyapunov stable. Then $f$ has a quasi horseshoe.}\\

{\bf Proof.} We may assume $h(f)=0$ since otherwise $f$ has a horseshoe.  By Theorem 2.1 there are $J(k_n,n)\in\mathcal J(\widetilde\omega)$ such that $J(k_{n+1},n+1)\subset J(k_n,n)$ and $\bigcap _{n\ge 0}J(k_n,n)=J$ is a non degenerate wandering interval; otherwise $f|_{\widetilde\omega}$ would be Lyapunov stable. Let $U_n:=K(k_{2n},2n)$  and $V_n:=K(k_{2n+1},2n+1)$, $n\in\mathbb N_0$. By Lemma 2.2 there are numbers $m_k$  such that $U_k, V_k, m_k$, $k\in\mathbb N_0$,  form a quasi horseshoe for $f$ with $\varepsilon = |J|$, the length of $J$.
$\hfill\Box$\\

{\bf Theorem 3.4.} {\it Let $f, f_k\in\mathcal C$, $k\in\mathbb N$, be surjective maps such that $f_{1,\infty}$ converges uniformly to $f$. If $f$ has an infinite $\omega$-limit set with isolated points then $(I, f_{1,\infty})$ is LYC.}\\

{\bf Proof.} Let $p\in\widetilde\omega:=\omega _f(z)$ be an isolated point, and let $J_p$ with endpoints $p$ and $a_p$ be as in Lemma 2.4. We show that there are sequences of compact intervals $K_j, P_j$, and positive integers $r_j, q_j$ such that
\begin{equation}
\label{nova420}
p\in P_j, \ f^{r_j}(K_j) \ \text{is a neighborhood of} \ K_j\cup P_{j}, \ f^{q_j}(P_j) \ \text{is a neighborhood of} \ K_{j+1}, \ \text{and}  \  r_j|q_j, \ j\in\mathbb N.
\end{equation}
To see this put $K_1=K(k_1,1)$. Let $J(k_n,n)$ be the intervals containing $p$ so that $\bigcap _{n\ge 0} J(k_n,n)=:J_p$. Since $\widetilde\omega\subset J(k_0,0)$,  Lemma 2.2  implies $f^{r_1}(K_1)\supset \widetilde\omega$ (where $r_1=4$), and since $J_p$ is a subset of the convex hull of $\widetilde\omega$, $f^{r_1}(K_1)$  must contain  infinitely many points from the trajectory of $z$ hence  a small neighborhood of $p$; denote it $P_1$. By Lemma 2.4  get $q_1$ divisible by $r_1$ such that $f^{q_1}(P_1)\supset J(k_0,n_2)\supset J_p$, and put $K_2:=K(0,n_2)$. By induction we get (\ref{nova420}) such that  $r_j$ are powers of $2$. Denote 
\begin{equation}
\label{nova421}
B_j:=K_j, f(K_j), f^2(K_j),\cdots ,f^{r_j-1}(K_j),\ \text {and} \  D_j:=P_j, f(P_j),f^2(P_j), \cdots ,f^{q_j-1}(P_j),  \ j\in\mathbb N.
\end{equation}
and consider the itinerary
\begin{equation}
\label{405}
\underbrace {B_1,B_1,\cdots , B_1}_\text{$c_1${ -times}},X_1, \underbrace {B_2,B_2,\cdots , B_2}_\text{$c_2${ -times}},X_2, \cdots , \underbrace {B_k,B_k,\cdots , B_k}_\text{$c_k${ -times}},X_k, \cdots, 
\end{equation}
where $c_k\in\mathbb N$, and 
\begin{equation}
\label{406}
X_k=D_k=:X^0_k \ \ \text{or} \ \ X_k=\underbrace{B_k,B_k, \cdots ,B_k}_\text{$q_k/r_k$ -times}=:X^1_k,
\end{equation}
so that the blocks  $X^0_k$ and $X^1_k$ have the same length $q_k$. By (\ref{nova420}) the above condition is correct. Let $\Sigma^\prime _2\subset\{0,1\}^\mathbb N$ be an uncountable set such that any two distinct sequences from $\Sigma _2^\prime$ have different coordinates at infinitely many places.  For $\beta =\{b_k\}_{k\ge 1}$ in $\Sigma _2^\prime$ let $x_\beta$ be a point in $I$ with trajectory (\ref{405}) such that $X_k=X_k^{b_k}$, $k\in\mathbb N$. If the numbers $a_k$ increase sufficiently rapidly then (\ref{405}) is the itinerary of a point $x^\prime_\beta$ for the nonautonomous system $f_{1,\infty}$, similarly as in the proof of Theorem A. Then
$S=\{x^\prime_\beta; \beta\in \Sigma _2^\prime\}$ is an uncountable scrambled set hence $(I,f_{1,\infty})$ is LYC, with $\varepsilon=|J_p|$.
$\hfill\Box$
\medskip

{\bf Proof of Theorem B.}  The result follows by Theorems 3.3 and 3.2 since $f$ has a horseshoe if $h(f)>0$.

$\hfill\Box$
\medskip

{\bf Proof of Theorem C.} By Theorem B we may assume that $h(f)=0$. Since $f$ is LYC, it has an infinite $\omega$-limit set $\widetilde\omega$ such that $f$ is not Lyapunov stable on it, see \cite{FShSm}.   If $\widetilde\omega$ has isolated points then the result follows by Theorem 3.4. Otherwise $\widetilde\omega$ is a minimal set;  apply Theorem  B. 
$\hfill\Box$

\section{Concluding remarks}

There are open problems related to our results.  We point out two of them. We assume $(I,f_{1,\infty})$ is a surjective system converging uniformly to $(I,f)$.
Then $f$ can be the identity map even if $(I,f_{1,\infty})$ is  chaotic, see, e.g., \cite{FPS}. In \cite{C} it is proved that if $(I,f_{1,\infty})$ is LYC then $f$ is LYC provided it has the shadowing property. But this condition eliminates maps $f$ with $h(f)=0$, see \cite{GK}. On the other hand, by Theorem B, if $h(f)>0$, then $f_{1,\infty}$ must be even DC1. 
\smallskip

{\it Problem 1.} Assume $(I,f_{1,\infty})$ is LYC and $h(f_{1,\infty})=0$. Find a condition for $f_{1,\infty}$  that is necessary and sufficient  for $f$ to be LYC.
\smallskip

 Uniform convergence of  $f_{1,\infty}$ to a map in $\mathcal C$ is essential to ensure that $h(f_{1,\infty})>0$ implies $(I,f_{1,\infty})$ is DC1: in \cite{SmSt} there is an example of a skew-product map $F:I^2\to I^2$ with $h(F)>0$ which is DC2, but not DC1. Recently T. Downarowicz \cite{D} proved that $h(f)>0$ implies DC2, for every $f\in\mathcal C(X)$. Recall that $(X,f)$ is DC2 if there is an uncountable set $S$ such that, for every distinct $x, y\in S$, $\Phi_{xy}<\Phi_{xy}^*\equiv 1$, cf. (\ref{equ12}).
 \smallskip

{\it Problem 2.} Assume $(I,f_{1,\infty})$ has positive topological entropy and $f_{1,\infty}$  converges {\it pointwise}  to a map in $\mathcal C$. Is it DC2? We conjecture that  $(I,f_{1,\infty})$ must have a DC2-pair.

\bigskip
{\bf Acknowledgement.} The author would like to thank prof. J. Sm\'{\i}tal for fruitful discussions and valuable comments.


\begin{thebibliography}{10}





\bibitem{BS} A. M. Bruckner and J. Sm\'{\i}tal, {\it A characterization of $\omega$-limit sets of maps of the interval with zero topological entropy\/}, Ergodic Theory \& Dynam. Systems {\bf13} (1993), 7--19. MR1213076 (94k:26006)

\bibitem{C}  J. C\'anovas, {\it Li-Yorke chaos in a class of nonautonomous discrete systems\/}, J. Difference Equ. Appl.  {\bf 17} (2011), 479--486. MR2783362 (2012c:37032)

\bibitem{D} T. Downarowicz, {\it Positive entropy implies distributional chaos DC2\/}, Proc. Amer. Math. Soc. 142 (2014),137--149. MR3119189

\bibitem{Dv} J. Dvo\v r\'akov\'a, {\it Chaos in nonautonomous discrete dynamical systems\/},  Commun. Nonlin. Sci. Numer. Simulat. {\bf 17} (2012), 4649-4652. MR2960260

\bibitem{FShSm} V. V. Fedorenko, A. N. \v Sarkovskii and J. Sm\'{\i}tal, {\it Characterizations of weakly chaotic maps of the interval\/}, Proc. Amer- Math. Soc. 110 (1990), 141--148. MR1017846 (91a:58148)

\bibitem{FPS} G.-L. Forti, L. Paganoni, and J. Sm\'{\i}tal, {\it Dynamics of homeomorphisms on minimal sets generated by triangular mappings\/}.  Bull. Austral. Math. Soc. 59 (1999), 1--20.  MR1672771 (99m:54029) 


\bibitem{KS} S. Kolyada and \v L. Snoha,  {\it Topological entropy of nonautonomous dynamical systems\/}, 
Random \& Comput.  Dynamics {\bf 4} (1996),  205--233. MR1402417 (98f:58126)

\bibitem{GK}  M. Kuchta,  {\it Shadowing property of continuous maps with zero topological entropy\/}, Proc. Amer. Math. Soc. {\bf 119} (1993), 641--648. MR1165058 (93k:58127)

\bibitem{OS} P. Oprocha and M. \v Stef\'ankov\'a, {\it Specification property and distributional chaos almost everywhere,}  Proc. Amer. Math. Soc., {\bf 136} (2008), 3931 - 3940. MR2425733 (2009i:37023)


\bibitem{ScSm} B. Schweizer and J. Sm\'{\i}tal, {\it Measures of chaos and a spectral decomposition of dynamical systems on the interval\/}, Trans. Amer. Math. Soc.  {\bf 344} (1994), 737--754. MR1225094 (94k:58091)

\bibitem{Sh} A. N. \v Sarkovskii, {\it Attracting sets containing no cycles,} Ukrain. Mat. \v Z. {\bf 20} (1968), no. 1, 136 - 142. MR0225314

\bibitem{S} J. Sm\'{\i}tal, {\it Chaotic functions with zero topological entropy\/}, Trans. Amer. Math. Soc. {\bf 297} (1986), 269--282. MR0849479 (87m:58107)

\bibitem{SmSt} J. Sm\'{\i}tal  and M. \v {S}tef\'ankov\'a, {\it Distributional chaos for triangular maps\/}, Chaos, Solitons and Fractals {\bf 21}, (2004), 1125--1128. MR2047330 (2005a:37017)


\end{thebibliography}
\end{document}